\documentclass[11pt]{amsart}
\usepackage{geometry}                % See geometry.pdf to learn the layout options. There are lots.
\usepackage{graphicx}
\usepackage{amssymb}
\usepackage{epstopdf}
\usepackage{breqn}
\usepackage{cite}
\title{An Analysis of the Replicator Dynamics for an Asymmetric Hawk-Dove Game}
\author{Ikjyot Singh Kohli}
\author{Michael C. Haslam}
\address{Department of Mathematics and Statistics, York University, Toronto, Ontario}
\address{Department of Mathematics and Statistics, York University, Toronto, Ontario}
\email{isk@mathstat.yorku.ca}
\email{mchaslam@mathstat.yorku.ca}
\date{March 16, 2017}                                           % Activate to display a given date or no date

\begin{document}

\begin{abstract}
In this paper, we analyze using a dynamical systems approach, the replicator dynamics for the asymmetric Hawk-Dove game in which there are a set of four pure strategies with arbitrary payoffs. We give a full account of the equilibrium points, their stability, and derive the Nash equilibria. We also give a detailed account of the local bifurcations that the system exhibits based on choices of the typical Hawk-Dove parameters $v,c$.  We also give details on the connections between the results found in this work and those of the standard two-strategy Hawk-Dove game. We conclude the paper with some examples of numerical simulations that further illustrate some global behaviour of the system. 
\end{abstract}

\maketitle

\section{Introduction}
The Hawk-Dove game is one of the first examples of a pairwise game that was used to model the conflict between animals \cite{smith1982evolution}. The basic idea is that ``Hawks'' and ``Doves'' represent two types of \emph{behaviours} (actions or pure strategies) that could be exhibited by animals of the same species \cite{webb1}. In the standard Hawk-Dove game, individuals can use one of two possible pure strategies. In one case, they can be aggressive/a ``Hawk'', which is typically denoted by $H$, or be non-aggressive/a ``Dove'', which is typically denoted by $D$. Then, at various times, individuals in this population can have a conflict over a resource which has value $v$, where the winner of the conflict gets the resource, which the loser pays a cost $c$. 

The hawk-dove game has been studied in the context of replicator dynamics a number of times over the past number of years. Some examples of these studies include \cite{bomze1983lotka, roca2009evolutionary, taylor2004evolutionary, ficici2000game, samuelson2002evolution, friedman1998economic, fudenberg2006evolutionary, johnstone2001eavesdropping, dall2004behavioural, hofbauer1998evolutionary, wydick2007games, ahmed2006dynamical, doebeli2005models, traulsen2007pairwise, weibull1997evolutionary}.

In replicator dynamics, it is assumed that individuals are programmed to use only pure strategies from a finite set $\mathbf{S} = \{s_{1}, \ldots, s_{k}\}$. It can be shown \cite{webb1} that the dynamical evolution of the proportion of individuals using strategy $s_{i}$, $x_i$, is given by:
\begin{equation}
\label{eq:dyn1}
\dot{x}_{i} = \left[\pi(s_i, \mathbf{x}) - \bar{\pi}(\mathbf{x})\right]x_{i},
\end{equation}
where $\pi(s_i, \mathbf{x})$ is the payoff to individuals using strategy $s_i$, while $\bar{\pi}(\mathbf{x})$ is known as the average payoff and is defined as
\begin{equation}
\label{eq:avpayoff1}
\bar{\pi}(\mathbf{x}) = \sum_{i=1}^{k} x_{i} \pi(s_i, \mathbf{x}).
\end{equation}
Further to Eq. \eqref{eq:dyn1}, one also has the constraint:
\begin{equation}
\label{eq:constr1}
\sum_{i=1}^{k} x_{i} = 1.
\end{equation}

In this paper, we wish to consider an asymmetric pairwise Dove-Hawk game. Following \cite{webb1}, specifically, this is where two individuals are contesting ownership of a territory that one of them controls. One assumes that the value of the territory and costs of contest are the same for both players. Unlike the standard Hawk-Dove game described above , players can now condition their behaviour on the role that they occupy, which is typically denoted as \emph{owner} or \emph{intruder}. So, the pure strategies now take the form \emph{play Hawk if owner}, and \emph{play Dove if intruder}, which we will denote by $HD$. Therefore, there are a set of four pure strategies:
\begin{equation}
\label{eq:strat4}
\mathbf{S} = \left\{HH,HD,DH,DD\right\}.
\end{equation}

From these arguments, it can be shown \cite{webb1}, that the payoff matrix is given by Table \ref{table1}.
\begin{table}[h]
\begin{center}\begin{tabular}{|c|c|c|c|c|}\hline  & $HH$ & $HD$ & $DH$ & $DD$ \\\hline $HH$ & $(v - c)/2, (v - c)/2$ & $(3 v - c)/4, (v - c)/4$ & $(3 v - c)/4, (v - c)/4$ & $v, 0$ \\\hline $HD$ & $(v - c)/4$, $(3 v - c)/4$ & $v/2, v/2$ & $(2 v - c)/4, (2 v - c)/4$ & $3 v/4, v/4$ \\\hline $DH$ & $(v - c)/4, (3 v - c)/4$ & $(2 v - c)/4, (2 v - c)/4$ & $v/2, v/2$ & $3 v/4, v/4$ \\\hline $DD$ & $0, v$ & $v/4, 3 v/4$ & $v/4, 3 v/4$ & $v/2,v/2$ \\\hline \end{tabular} \caption{The payoff matrix for the asymmetric Hawk-Dove game.}
\end{center}
\label{table1}
\end{table}

We note that the replicator dynamics of this four-strategy asymmetric hawk-dove game has not been analyzed from a dynamical systems perspective in the literature to the best of the authors' knowledge. However, some examples of related asymmetric games can be found in  \cite{Mesterton-Gibbons1992, matsumura1998game, nakamaru2003can, he2014evolutionary, mcavoy2015asymmetric, uehara2010global, sekiguchi2015fixation,benz2005game,cressman2003evolutionary}.

\section{The Dynamical Equations}
Let us denote by $(x,y,z,w)$ the proportion of individuals who use strategies $HH, HD, DH$ and $DD$ respectively. Then,from the payoff matrix in Table \ref{table1} and Eqs. \eqref{eq:dyn1}-\eqref{eq:avpayoff1}, the replicator dynamics are given by the following dynamical system:
\begin{eqnarray}
\label{eq:xdot}
\dot{x} &=& \frac{1}{4} x \left[-c (2 x+y+z)-4 \bar{\pi}(\mathbf{x})+v (4 w+2 x+3 (y+z))\right], \\
\label{eq:ydot}
\dot{y} &=& \frac{1}{4} y \left[-c (x+z)-4  \bar{\pi}(\mathbf{x})+v (3 w+x+2 (y+z))\right], \\
\label{eq:zdot}
\dot{z} &=& \frac{1}{4} z \left[-c (x+y)-4  \bar{\pi}(\mathbf{x})+v (3 w+x+2 (y+z))\right], \\
\label{eq:wdot}
\dot{w} &=& \frac{1}{4} w \left[v (2 w+y+z)-4  \bar{\pi}(\mathbf{x})\right], \\
\end{eqnarray}
where
\begin{equation}
 \bar{\pi}(\mathbf{x}) = \frac{1}{2} \left[v -c (x+y) (x+z)\right],
\end{equation}
and from Eq. \eqref{eq:constr1},
\begin{equation}
\label{eq:constr2}
x + y + z + w = 1.
\end{equation}

This four-dimensional dynamical system can be reduced to three dimensions if we set via Eq. \eqref{eq:constr2}, $w = 1 - x - y - z$. Therefore, in what follows, we will study the following \emph{unconstrained} three-dimensional system:
\begin{eqnarray}
\label{eq:xdot2}
\dot{x} &=& \frac{1}{4} x \left[c \left(2 x^2+2 x (y+z-1)+2 y z-y-z\right)-v (2 x+y+z-2)\right], \\
\label{eq:ydot2}
\dot{y} &=& -\frac{1}{4} y \left[v (2 x+y+z-1)-c (2 x+2 y-1) (x+z)\right], \\
\label{eq:zdot2}
\dot{z} &=& -\frac{1}{4} z \left[v (2 x+y+z-1)-c (x+y) (2 x+2 z-1)\right].
\end{eqnarray}

\section{A Local Stability Analysis}
From Eqs. \eqref{eq:xdot2}-\eqref{eq:zdot2}, we now present the equilibrium points along with their eigenvalues and local stability. The Jacobian matrix, denoted $J_{ij}$ corresponding to this dynamical system is a $3\times3$ matrix, whose entries we list as follows:
\begin{eqnarray}
J_{11} &=& \frac{1}{4} \left(c \left(6 x^2+4 x (y+z-1)+2 y z-y-z\right)-v (4 x+y+z-2)\right), \\
J_{12} &=& \frac{1}{4} x (c (2 x+2 z-1)-v), \\
J_{13} &=& \frac{1}{4} x (c (2 x+2 y-1)-v), \\
J_{21} &=& \frac{1}{4} y (c (4 x+2 y+2 z-1)-2 v), \\
J_{22} &=& \frac{1}{4} (c (2 x+4 y-1) (x+z)-v (2 x+2 y+z-1)), \\
J_{23} &=& \frac{1}{4} y (c (2 x+2 y-1)-v), \\
J_{31} &=& \frac{1}{4} z (c (4 x+2 y+2 z-1)-2 v), \\
J_{32} &=& \frac{1}{4} z (c (2 x+2 z-1)-v), \\
J_{33} &=& \frac{1}{4} (c (x+y) (2 x+4 z-1)-v (2 x+y+2 z-1)).
\end{eqnarray}

\subsection{Equilibrium Point 1}
The first equilibrium point was found to be 
\begin{equation}
\label{eq:P1}
P_{1}: (x,y,z) = (0,0,1) 
\end{equation}
The corresponding eigenvalues of $J_{ij}$ were found to be:
\begin{equation}
\label{eq:eig1}
\{\lambda_{1},\lambda_{2}, \lambda_{3}\} = \left\{-\frac{v}{4},-\frac{c}{4},\frac{v-c}{4}\right\}. 
\end{equation}
This point is a stable node if:
\begin{equation}
\{v > 0\} \cap \{c > v\}.
\end{equation}
It is an unstable node if:
\begin{equation}
\{v < 0\} \cap \{c < v\}.
\end{equation}
It is a saddle if
\begin{equation}
\left\{ \{v <0\} \cap \{v < c < 0\} \right\} \bigcup \left\{ \{v < 0\} \cap \{c > 0\}\right\} \bigcup \left\{ \{v > 0\} \cap \{c < 0\} \right\} \bigcup \left\{ \{v > 0\} \cap \{ 0 < c <v\}     \right\}.
\end{equation}

\subsection{Equilibrium Point 2}
The second equilibrium point was found to be:
\begin{equation}
\label{eq:P2}
P_{2}: (x,y,z)= \left(0, \frac{1}{2}, \frac{1}{2}\right).
\end{equation}
The corresponding eigenvalues of $J_{ij}$ were found to be:
\begin{equation}
\label{eq:eig2}
\{\lambda_{1},\lambda_{2}, \lambda_{3}\} = \left\{\frac{c}{8},\frac{1}{8} (c-2 v),\frac{1}{8} (2 v-c)\right\}.
\end{equation}
This point is neither a stable or unstable node. However, it is a saddle under the following conditions:
\begin{equation}
\left\{(v\leq 0\cap c<2 v)\cup (v>0\cap c<0)\right\} \bigcup \left\{v>0\cap0<c<2 v    \right\} \bigcup \left\{ v<0\cap 2 v<c<0  \right\}.
\end{equation}

\subsection{Equilibrium Point 3}
The third equilibrium point was found to be:
\begin{equation}
\label{eq:P3}
P_{3}: (x,y,z) = \left(0, \frac{v}{c}, \frac{v}{c}   \right).
\end{equation}
The corresponding eigenvalues of $J_{ij}$ were found to be:
\begin{equation}
\label{eq:eig3}
\{\lambda_{1},\lambda_{2}, \lambda_{3}\} = \left\{\frac{v}{4},-\frac{v (c-2 v)}{4 c},0\right\}.
\end{equation}
The single zero eigenvalue indicates that this equilibrium point is \emph{normally} hyperbolic, and the local stability can be determined through the non-zero eigenvalues by the invariant manifold theorem \cite{ellis}. In particular, this point is a stable node if:
\begin{equation}
v<0 \cap 2 v<c<0.
\end{equation}
It is an unstable node if:
\begin{equation}
v>0\cap 0<c<2 v.
\end{equation}
It is a saddle point under the following conditions:
\begin{equation}
\left\{v<0\cap (c<2 v\cup c>0)    \right\}  \bigcup \left\{ v>0\cap (c<0\cup c>2 v)   \right\}.
\end{equation}

\subsection{Equilibrium Point 4}
The fourth equilibrium point was found to be:
\begin{equation}
\label{eq:P4}
P_{4}: \left(x,y,z\right) = \left(0,1,0\right).
\end{equation}
The corresponding eigenvalues of $J_{ij}$ were found to be:
\begin{equation}
\label{eq:eig4}
\{\lambda_{1},\lambda_{2}, \lambda_{3}\}  = \left\{-\frac{v}{4},-\frac{c}{4},\frac{v-c}{4}\right\}.
\end{equation}
This point is a stable node if
\begin{equation}
\{v > 0\} \cap \{c > v\}.
\end{equation}
It is an unstable node if
\begin{equation}
\{v < 0 \} \cap \{c < v\}.
\end{equation}
It is a saddle point under the following conditions:
\begin{equation}
\left\{v<0 \ \cap v <c<0    \right\} \bigcup \left\{ v<0\cap c>0    \right\} \bigcup \left\{ v>0\cap c<0   \right\} \bigcup \left\{   v>0\cap 0<c<v   \right\}.
\end{equation}

\subsection{Equilibrium Point 5}
The fifth equilibrium point was found to be:
\begin{equation}
\label{eq:P5}
P_{5}: \left(x,y,z\right) = \left(1,0,0\right).
\end{equation}
The corresponding eigenvalues of $J_{ij}$ were found to be:
\begin{equation}
\label{eq:eig5}
\{\lambda_{1},\lambda_{2}, \lambda_{3}\}  = \left\{\frac{c-v}{2},\frac{c-v}{4},\frac{c-v}{4}\right\}.
\end{equation}
This point is a stable node if
\begin{equation}
\{v \in \mathbb{R}\} \cap \{c < v\}.
\end{equation}
It is an unstable node if
\begin{equation}
\{v \in \mathbb{R}\} \cap \{c > v\}.
\end{equation}
From Eq. \eqref{eq:eig5}, it can be seen that $\mathcal{P}_{5}$ is in fact never a saddle point of the dynamical system.

\subsection{Equilibrium Point 6}
The sixth equilibrium point was found to be:
\begin{equation}
\label{eq:P6}
P_{6}: \left(x,y,z\right) = \left( \frac{v}{c}, 0, 0\right).
\end{equation}
The corresponding eigenvalues of $J_{ij}$ were found to be:
\begin{equation}
\label{eq:eig6}
\{\lambda_{1},\lambda_{2}, \lambda_{3}\}  = \left\{0,0,\frac{v (v-c)}{2 c}\right\}.
\end{equation}
One sees that since $\lambda_{1} = \lambda_{2} = 0$, this point is manifestly non-hyperbolic. As such, its stability properties cannot be determined through the Jacobian matrix.

\subsection{Equilibrium Point 7}
The final equilibrium point was found to be:
\begin{equation}
\label{eq:P7}
P_{7}: \left(x,y,z\right) = \left(0,0,0\right).
\end{equation}
The corresponding eigenvalues of $J_{ij}$ were found to be:
\begin{equation}
\label{eq:eig7}
\{\lambda_{1},\lambda_{2}, \lambda_{3}\}  = \left\{\frac{v}{2},\frac{v}{4},\frac{v}{4}\right\}.
\end{equation}
This point is a stable node if
\begin{equation}
v  < 0.
\end{equation}
It is an unstable node if
\begin{equation}
v > 0.
\end{equation}
Further, this point is never a saddle point as can be seen from Eq. \eqref{eq:eig7}.

%One can also see that there exists a finite heteroclinic sequence from this analysis of the equilibrium points. For any $c > 0, v > c$, one has the following heteroclinic sequence:
%\begin{equation}
%P_{3}:(0,v/c,v/c) \rightarrow P_{2}:(0,1/2,1/2) \rightarrow P_{5}: (1,0,0),
%\end{equation}
%where in this case, $P_{3}$ is an unstable node, $P_{2}$ is a saddle, and $P_{5}$ is a stable node.

\section{Local Bifurcations}
With knowledge of the equilibrium points and their local stability as given in the previous sections, we now attempt to describe bifurcation behaviour exhibited by this dynamical system. Analyzing bifurcation behaviour is important as this determines the local changes in stability of the equilibrium points of the system. 

The mechanism for these bifurcations can be seen as follows. 

The linearized system in a neighbourhood of $P_{1}$ takes the form:
\begin{eqnarray}
\dot{x} &=& \frac{1}{4} x \left(v - c\right), \\
\dot{y} &=& -\frac{1}{4} y c, \\
\dot{z} &=& \frac{1}{4}x \left(c - 2v \right) + \frac{c-v}{4} y  - \frac{v}{4}z.
\end{eqnarray}
We see that $x$ destabilizes $P_{1}$ when $v = c$, $y$ destabilizes $P_{1}$ when $c = 0$, and that $z$ destabilizes $P_{1}$ when $v = c = 0$.

The linearized system in a neighbourhood of $P_{2}$ takes the form:
\begin{eqnarray}
\dot{x} &=& \frac{1}{8}x (2v - c), \\
\dot{y} &=& \frac{1}{8}x (c - 2v) + \frac{1}{8}y (c-v) - \frac{1}{8}vz, \\
\dot{z} &=& \frac{1}{8}x (c - 2v) -\frac{1}{8}yv + \frac{1}{8}(c-v)z.
\end{eqnarray}
We see that $x$ destabilizes $P_{2}$ along the line $c = 2v$, while $y$ and $z$ destabilize $P_{2}$ when $v = c = 0$.

The linearized system in a neighbourhood of $P_{3}$ takes the form:
\begin{eqnarray}
\dot{x} &=& 0, \\
\dot{y} &=& -\frac{v (c-2 v)}{4 c} x + \frac{v^2}{4 c} y + \frac{v (v-c)}{4 c} z, \\
\dot{z} &=&  -\frac{v (c-2 v)}{4 c} x + \frac{v (v-c)}{4 c} y + \frac{v^2}{4 c} z.
\end{eqnarray}
Therefore, $P_{3}$ is destabilized by $y$ and $z$ $\forall \ c \neq 0, v = 0$.

The linearized system in a neighbourhood of $P_{4}$ takes the form:
\begin{eqnarray}
\dot{x} &=& \frac{1}{4}x(v-c), \\
\dot{y} &=& \frac{1}{4}(c - 2v)x - \frac{v}{4}y + \frac{c-v}{4}z, \\
\dot{z} &=& -\frac{c}{4} z.
\end{eqnarray}
Therefore $x$ destabilizes $P_{4}$ along the line $v = c$. Further, $y$ destabilizes $P_{4}$ when $v = c = 0$. Finally, $z$ destabilizes $P_{4}$ when $c = 0$ for $v \in \mathbb{R}$.

The linearized system in a neighbourhood of $P_{5}$ takes the form:
\begin{eqnarray}
\dot{x} &=& \frac{c-v}{2}x  + \frac{c-v}{4}y + \frac{c-v}{4}z, \\
\dot{y} &=& \frac{c-v}{4} y, \\
\dot{z} &=& \frac{c-v}{4} z.
\end{eqnarray}
Therefore $P_{5}$ is destabilized by $x$, $y$, and $z$ along the line $c = v$.

The linearized system in a neighbourhood of $P_{6}$ takes the form:
\begin{eqnarray}
\dot{x} &=& \frac{v (v-c)}{2 c} x + \frac{v (v-c)}{4 c} y + \frac{v (v-c)}{4 c}, \\
\dot{y} &=& 0, \\
\dot{z} &=& 0.
\end{eqnarray}
Therefore, $x$ destabilizes $P_{6}$ whenever $v = 0$, or whenever $v = c$ (for $c \neq 0$).

The linearized system in a neighbourhood of $P_{7}$ takes the form:
\begin{eqnarray}
\dot{x} &=& \frac{v}{2} x, \\
\dot{y} &=& \frac{v}{4} y, \\
\dot{z} &=& \frac{v}{4} z.
\end{eqnarray}
We see that therefore, $P_{7}$ is destabilized by $x$, $y$, and $z$ whenever $v = 0$, for $c \in \mathbb{R}$.
% v =c 
% c = 0
% v = c = 0
% c = 2v
% v = c = 0 
% v = 0, for any c
From these calculations, we can therefore see that along $v = c$, for as one goes from $c  < 0$ to $c  > 0$, $P_{1}$ and $P_{4}$ go from being unstable nodes to a stable ones, and vice-versa, while $P_{3}$ goes from being a stable node to an unstable one. Whenever $c = 0, v \in \mathbb{R}$, and as one goes from $v < 0$ to $v > 0$, $P_5$ goes from being an unstable node to a stable node, while $P_7$ goes from being a stable node to an unstable one. Along the line $c = 2v$, as we go from $v < 0$ to $v > 0$, $P_{1}$ and $P_{4}$ go from being unstable nodes to stable nodes, while $P_5$ and $P_7$ go from being stable nodes to unstable nodes. Finally, whenever $v = 0, c \in \mathbb{R}$, as we go from $c < 0$ to $c > 0$, $P_{5}$ goes from being a stable node to an unstable one, while $P_1$ and $P_4$ go from being unstable nodes to stable ones.

\section{Nash Equilibria}
Determining the future asymptotic behaviour of the replicator dynamics is of importance since by Theorem 9.15 in \cite{webb1}, if $\mathbf{x}^{*}$ is an asymptotically stable fixed point of the dynamical system, then the symmetric strategy pair $[\sigma^{*}, \sigma^{*}] = [\mathbf{x}^{*}, \mathbf{x}^{*}]$ is a Nash equilibrium.   

Following \cite{arnolddyn}, we note that first, by Lyapunov's theorem, if all eigenvalues of the linear part of a vector field $v$ at a singular point have negative real part, the singular point is asymptotically stable.

From our stability analysis of the various equilibrium points in the preceding sections, we therefore observe the following Nash equilibria of the replicator dynamics depending on the choices of $v$ and $c$:
\begin{enumerate}
\item $v > 0, c > v$ $\Rightarrow$  $P_{1}$ is asymptotically stable $\Rightarrow [(0,0,1),(0,0,1)]$ is a Nash equilibrium.
\item $ v > 0, c > v$ $\Rightarrow$ $P_{4}$ is asymptotically stable $\Rightarrow [(0,1,0), (0,1,0)]$ is a Nash equilibrium. 
\item $v \in \mathbb{R}, c < v$ $\Rightarrow$ $P_{5}$ is asymptotically stable $\Rightarrow [(1,0,0), (1,0,0)]$ is a Nash equilibrium.
\item $c \in \mathbb{R}, v  < 0$ $\Rightarrow$ $P_{7}$ is asymptotically stable $\Rightarrow [(0,0,0), (0,0,0)]$ is a Nash equilibrium. 
\end{enumerate}

The existence of these Nash equilibria show that this asymmetric Hawk-Dove game produces rational behaviour in a population composed of players that are not required to make \emph{consciously} rational decisions. In other words, the population is stable when, given what everyone else is doing, no individual would get a better result by adopting a different strategy. This is the so-called \emph{population view} of a Nash equilibrium, which Nash himself described as the mass action interpretation \cite{webb1, nashphd}.

\section{Connections with The Two-Strategy Hawk-Dove Game}
%added this section based on referee report...
It is perhaps of interest to discuss our results found above in connection with the standard two-strategy hawk-dove game. Following \cite{webb1}, we note that the payoff matrix for such a game is given by
\begin{table}[h]
\begin{center}
\begin{tabular}{|c|c|c|}\hline  & $H$ & $D$ \\\hline $H$ & $(v-c)/2,  (v-c)/2$ & $v,0$ \\\hline $D$ & $0,v$ & $v/2,v/2$ \\\hline 
\end{tabular} 
\caption{Payoff matrix for the standard two-strategy Hawk-Dove game.}
\label{table2}
\end{center}
\end{table}

In this case, the replicator dynamics are a simple consequence of Eqs. \eqref{eq:dyn1}-\eqref{eq:avpayoff1}. Namely, let $z$ denote the proportion of individuals in the population that use strategy $H$ in Table \ref{table2}. Then, the replicator dynamics is governed by the single ordinary differential equation:
\begin{equation}
\label{eq:case2}
\dot{z} = \frac{c}{2} z \left(1-z\right) \left(\frac{v}{c} - z\right).
\end{equation}
Clearly, Eq. \eqref{eq:case2} has equilibrium points $z = 0$, $z = 1$, and $z = v/c$. Let us denote by $f(z)$ the right-hand-side of Eq. \eqref{eq:case2}. Then,
\begin{equation}
\label{eq:jacobian1d}
f'(z) = \frac{1}{2} \left[v - 2vz + cz \left(-2 + 3z\right)\right].
\end{equation}
Clearly, when $z = 0$, $f'(z) = v/2$, which is negative when $v < 0$ and positive when $v > 0$. Therefore, the point $z=0$ is a stable node when $v<0$, and an unstable node when $v>0$. Further, when $z = 1$, $f'(z) = (c-v)/2$. In this case, the point $z=0$ is a stable node for $c \in \mathbb{R}$, and $v > c$. Further, it is an unstable node for $c \in \mathbb{R}$ and $v < c$. Finally, when $z = v/c$, we have that $f'(z) = v(v-c)/2c$. This point is a stable node when $v < 0$ and $v < c < 0$, or, when $v >0$ and $c < 0$ or $c > v$. It is unstable node when $v < 0$ and $c < v$ or $c > 0$, or when $v > 0$ and $0 < c < v$. 

Comparing these cases to the Nash equilibria we found in the full asymmetric game. We see that the case when $z = 0$ corresponds to the case of Equilibrium Point 7 above, where $[(0,0,0), (0,0,0)]$ was a Nash equilbrium. The case $z = v/c$ in this example corresponds to Equilibrium Points 1 and 4 above, where $[(0,0,1),(0,0,1)]$ and $[(0,1,0), (0,1,0)]$ were both found to be Nash equilibria of the full asymmetric replicator dynamics. Certainly, this shows that $z \to v/c$ for any initial population that is not at an equilibrium point.

%add some interpretations here....

%%Is equilibrium point 6 asymptotically stable??
%Following Proposition A1 in \cite{leblanc3}, if $\phi_{t}$ is a flow on $\mathbb{R}^{n}$ with $S$ an invariant set, and if $Z: S \to \mathbb{R}$ is a $C^{1}$ function whose range is the interval $(a,b)$, where $a \in \mathbb{R} \cup \{-\infty\}$, $b \in \mathbb{R} \cup \{+\infty\}$ and $a < b$, then if $Z$ is decreasing on orbits in $S$, for all $\mathbf{x} \in S$ we have that $\omega(\mathbf{x}) \subseteq \left\{\mathbf{s} \in \bar{S} \backslash S: \lim_{\mathbf{y} \to \mathbf{s}} Z (\mathbf{y}) \neq  b\right\}$,  $\alpha(\mathbf{x}) \subseteq \left\{\mathbf{s} \in \bar{S} \backslash S: \lim_{\mathbf{y} \to \mathbf{s}} Z (\mathbf{y}) \neq  a\right\}$. 
%
%Further, by Proposition 4.1 in \cite{ellis}, consider a DE $\mathbf{x}' = \mathbf{f(x)}$ on $\mathbb{R}^{n}$ with flow $\phi_{t}$. Let $Z: \mathbb{R}^{n} \to \mathbb{R}$ be a $C^{1}$ function with satisfies $Z' = \alpha Z$, where $\alpha: \mathbb{R}^{n} \to \mathbb{R}$ is a continuous function. Then the subsets of $\mathbb{R}^{n}$ defined by $Z > 0$, $Z = 0$, $Z < 0$ are invariant sets of the flow $\phi_{t}$.
%
%Applying this latter proposition to Eqs. \eqref{eq:xdot2}-\eqref{eq:zdot2}, we see that $x > 0$, $y > 0$, and $z > 0$, and $x = 0$, $y = 0$, $z = 0$ are invariant sets of the dynamical system. 

\section{Some Numerical Simulations}
In this section, we present some numerical simulations of the work above. These simulations were completed in MATLAB using the ODE23s solver with a variety of initial conditions which are denoted with asterisks in the plots that follow. 

\begin{figure}[h]
%\begin{center}
\includegraphics[scale=0.40]{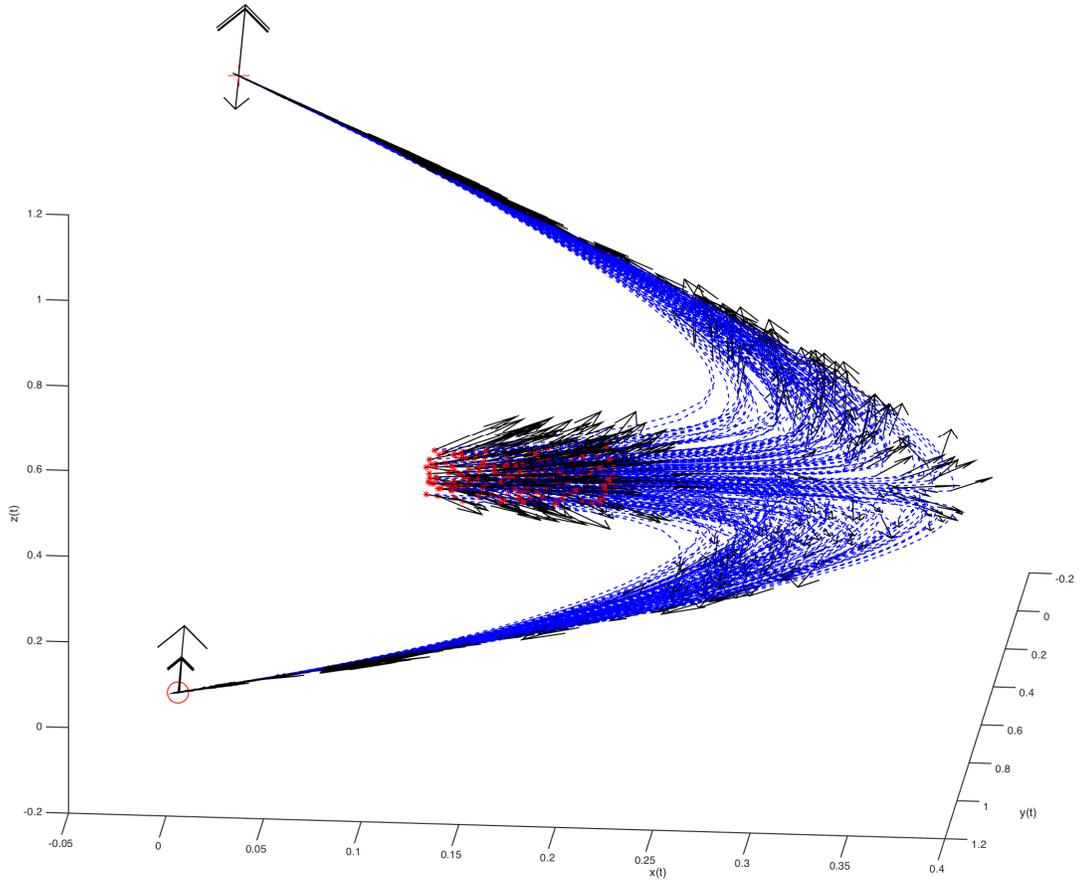}
\caption{Results where $v = 0.1, c = 0.2$. The red cross denotes the equilibrium point $P_{1}:(x,y,z) = (0,0,1)$, and the red circle denotes the equilibrium point $P_{4}:(x,y,z)=(0,1,0)$.}
\label{fig:fig1}
%\end{center}
\end{figure}

\begin{figure}[h]
%\begin{center}
\includegraphics[scale=0.40]{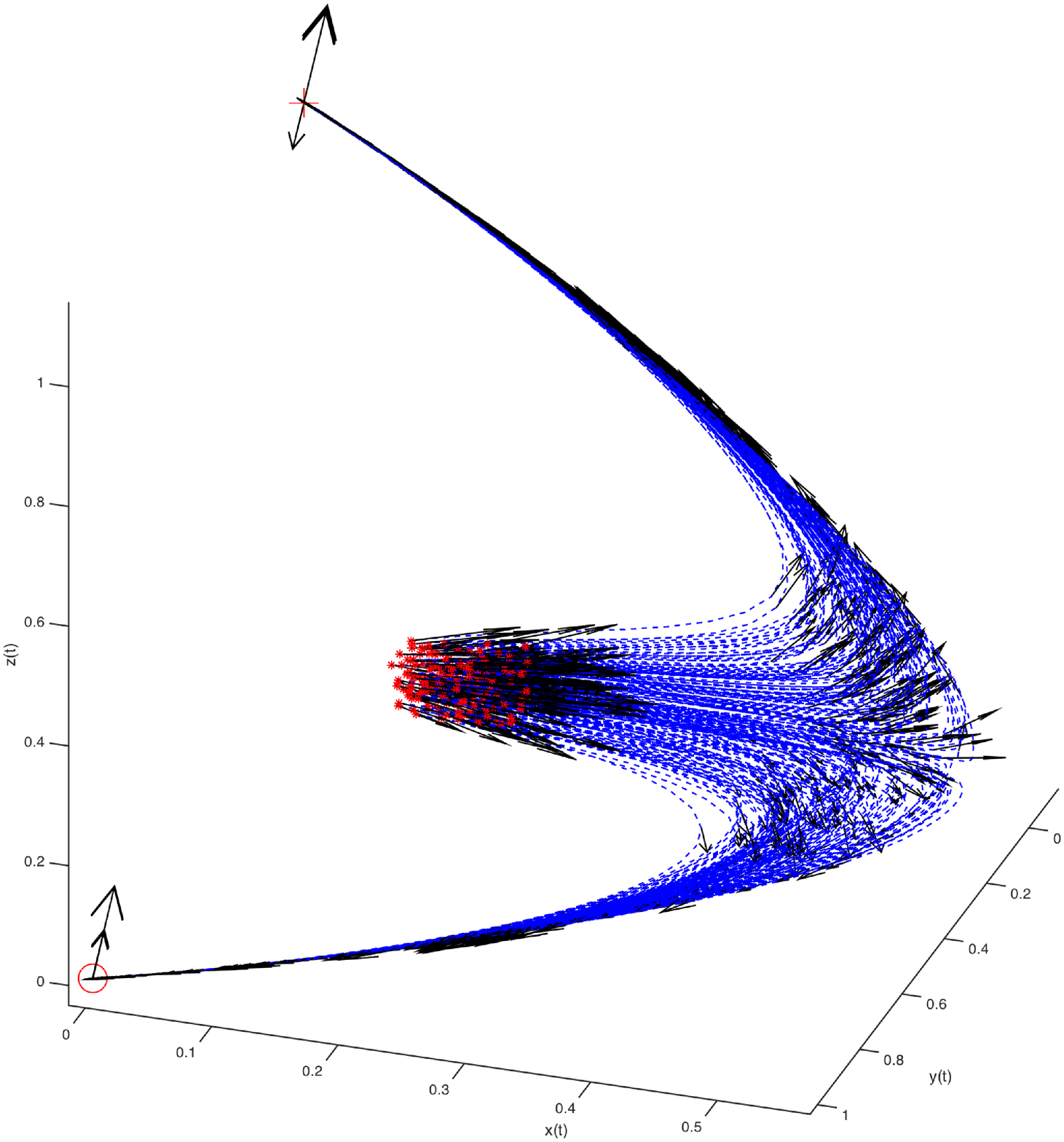}
\caption{Results where $v = 0.2, c = 0.3$. The red cross denotes the equilibrium point $P_{1}:(x,y,z) = (0,0,1)$, and the red circle denotes the equilibrium point $P_{4}:(x,y,z)=(0,1,0)$.}
\label{fig:fig2}
%\end{center}
\end{figure}

\begin{figure}[h]
%\begin{center}
\includegraphics[scale=0.40]{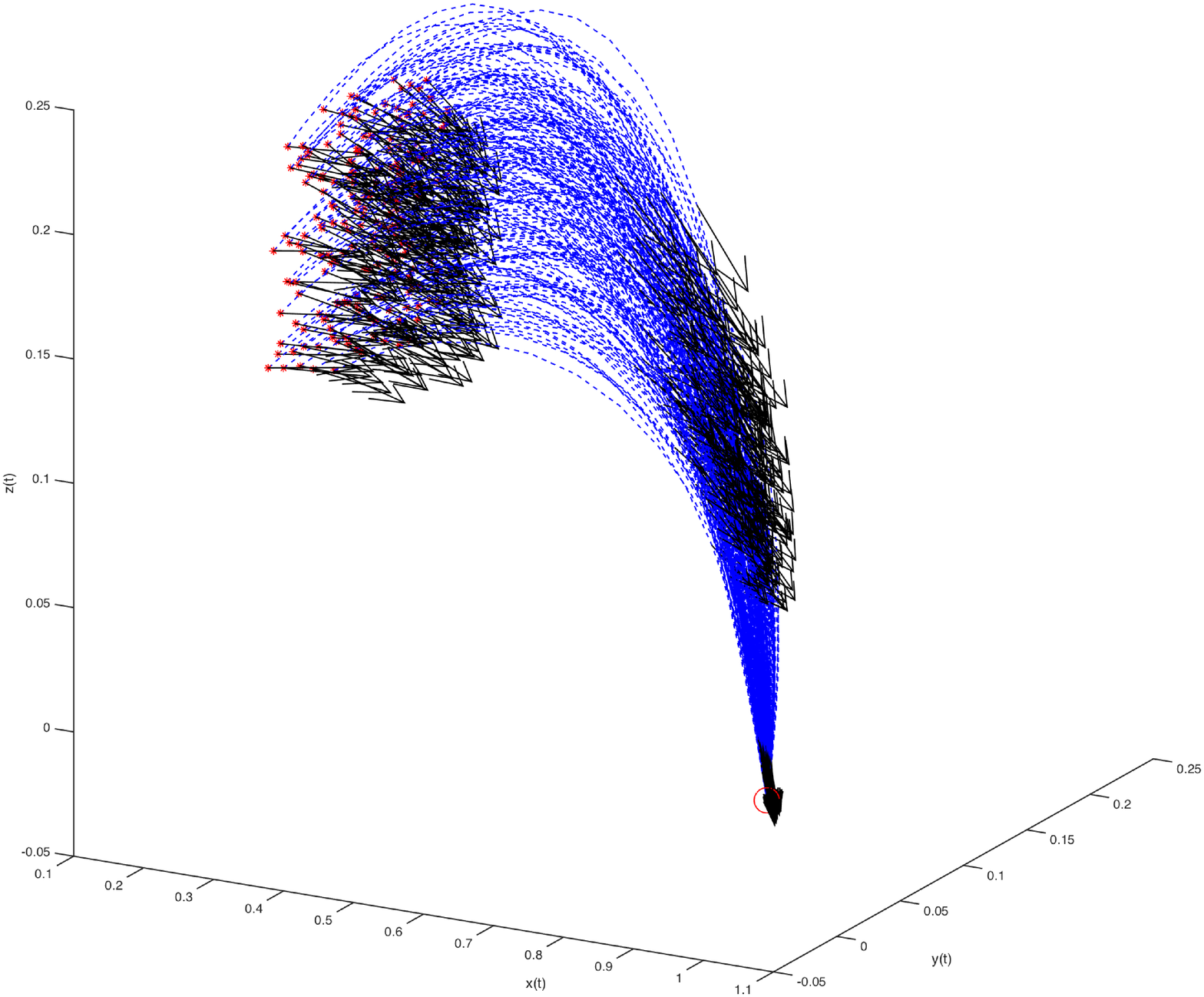}
\caption{Results where $v = 0.2, c = 0.1$. The red circle denotes the equilibrium point $P_{5}:(x,y,z) = (1,0,0)$. }
\label{fig:fig3}
%\end{center}
\end{figure}

\begin{figure}[h]
%\begin{center}
\includegraphics[scale=0.40]{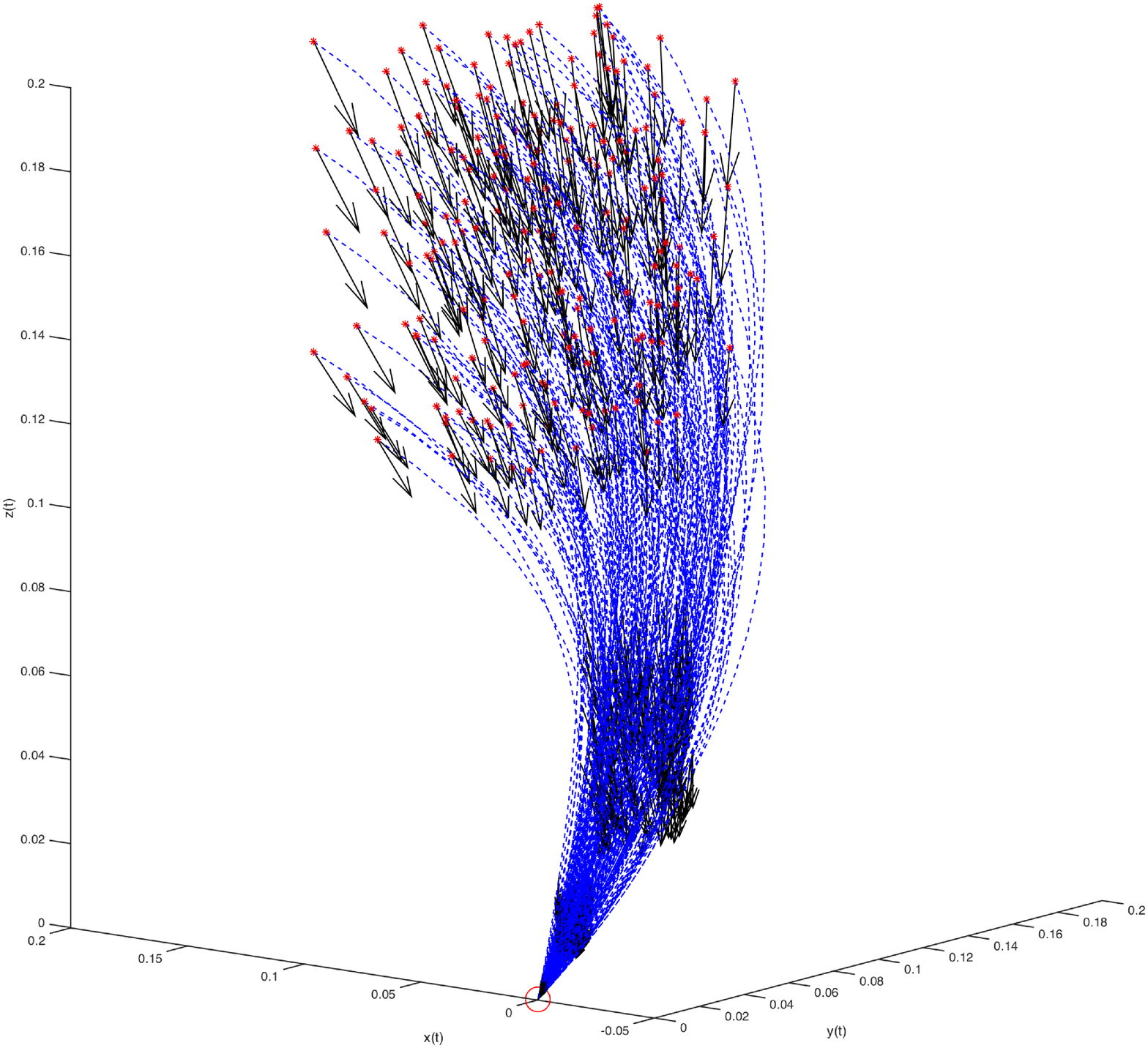}
\caption{Results where $v = -0.1, c = 0.2$. The red circle denotes the equilibrium point $P_{7}:(x,y,z) = (0,0,0)$.}
\label{fig:fig4}
%\end{center}
\end{figure}

\begin{figure}[h]
%\begin{center}
\includegraphics[scale=0.40]{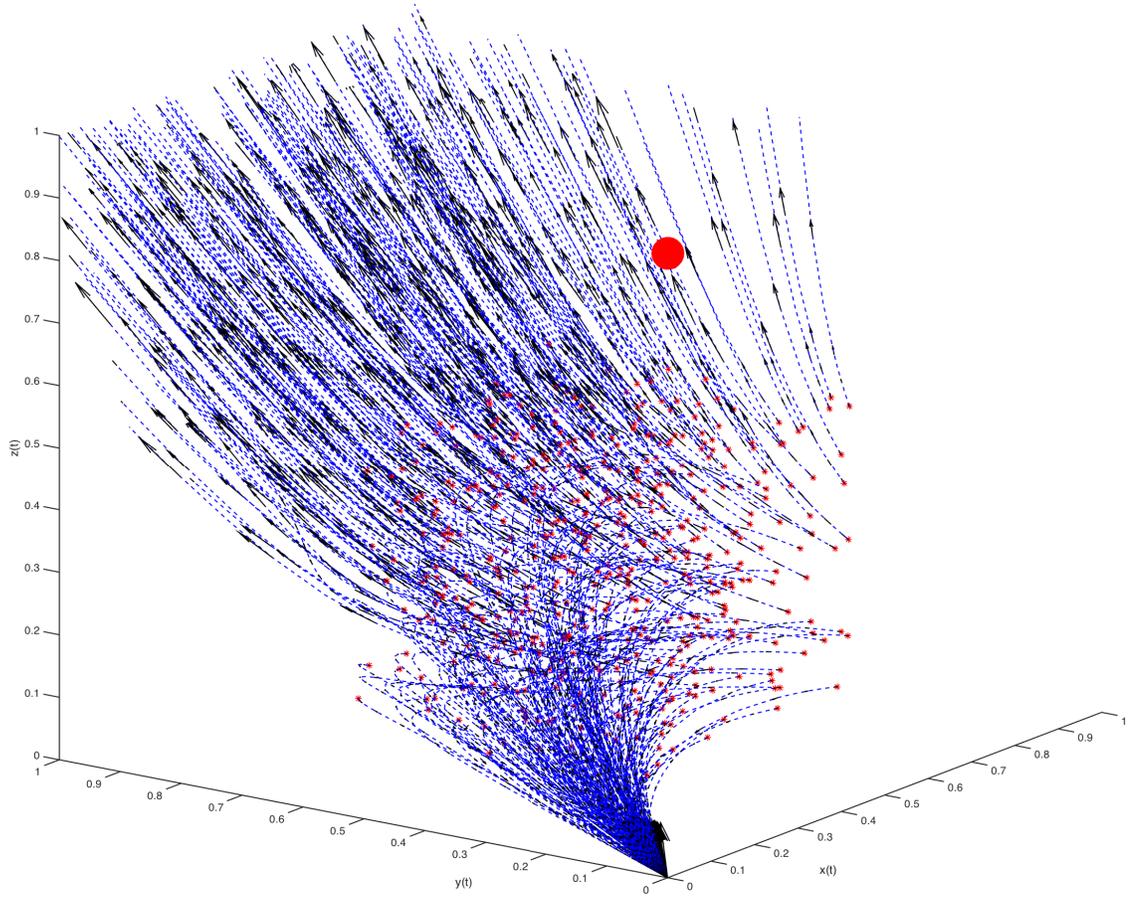}
\caption{Results where $v = -0.2, c = -0.1$. The shaded red circle denotes the equilibrium point $P_{1}:(x,y,z) = (0,0,1)$. One can see that indeed under these choices for $v$ and $c$, $P_1$ is indeed a saddle as predicted by our stability analysis.}
\label{fig:fig5}
%\end{center}
\end{figure}

In Fig. \ref{fig:fig1}, we assume that $v = 0.1, c = 0.2$, in Fig. \ref{fig:fig2}, we assume that $v = 0.2, c= 0.3$, in Fig. \ref{fig:fig3}, we assume that $v = 0.2, c = 0.1$, in Fig. \ref{fig:fig4}, we assume that $v = -0.1, c = 0.2$, and in Fig. \ref{fig:fig5}, $v = -0.2, c = -0.1$.
\newpage
\newpage
\newpage
\newpage
\newpage
\newpage

\section{Conclusions}
In this paper, we analyzed using a dynamical systems approach, the replicator dynamics for the asymmetric Hawk-Dove game in which there are a set of four pure strategies with arbitrary payoffs. We gave a full account of the equilibrium points, their stability, and derived the Nash equilibria. In particular, we found that if $v > 0, c>0$, then the strategy pairs $[HD,HD]$ and $[DH,DH]$ are Nash equilibria. If $ v \in \mathbb{R}, c < v$, then the strategy pair $[HH,HH]$ is a Nash equilibrium. Finally, if $c \in \mathbb{R}, v < 0$, then the strategy pair $[DD,DD]$ is a Nash equilibria. We also gave a detailed account of the local bifurcations that the system exhibits based on choices of the typical Hawk-Dove parameters $v,c$.  We also gave details on the connections between the results we found and those of the standard two-strategy Hawk-Dove game. We concluded the paper with some examples of numerical simulations that further illustrate some global behaviour of the system.

\section{Acknowledgements}
This research was partially supported by a grant given to MCH from the Natural Sciences and Engineering Research Council of Canada. The authors would also like to thank the anonymous referee for helpful comments and suggestions.
\bibliographystyle{apalike}
\bibliography{sources}

\end{document}